\documentclass[11pt]{article}
\usepackage[margin=1in]{geometry}
\usepackage{xcolor}
\usepackage{amsmath,amssymb,amsthm,mathtools}
\usepackage[T1]{fontenc}
\usepackage{lmodern,microtype,enumitem,longtable}
\usepackage[colorlinks=true,linkcolor=blue,citecolor=blue,urlcolor=blue]{hyperref}
\newtheorem{theorem}{Theorem}[section]
\newtheorem{lemma}[theorem]{Lemma}
\newtheorem{proposition}[theorem]{Proposition}
\newtheorem{corollary}[theorem]{Corollary}
\theoremstyle{remark}\newtheorem{remark}[theorem]{Remark}
\newcommand{\T}{\mathbb T}
\newcommand{\R}{\mathbb R}
\newcommand{\norm}[1]{\left\|#1\right\|}

\newcommand{\calM}{\mathcal M}
\newcommand{\calA}{\mathcal A}

\DeclareMathOperator{\spanop}{span}
\numberwithin{equation}{section}

  \chardef\forshowkeys=0
  \chardef\showllabel=0
  \chardef\refcheck=0
  \chardef\sketches=0
  \chardef\figures=1

\ifnum\forshowkeys=1

  \usepackage[notref,notcite,color]{showkeys}
\fi

\definecolor{mygray}{rgb}{.6,.6,.6}
\definecolor{myblue}{rgb}{0.9, 0, 1}
\definecolor{colorforkeys}{rgb}{1.0,0.0,0.0}

\definecolor{colorigor}{rgb}{.5, 0.2, 0.8}

\definecolor{colorgggg}{rgb}{0.1,0.5,0.3}
\definecolor{colorllll}{rgb}{0.0,0.7,0.0}
\definecolor{colorhhhh}{rgb}{0.3,0.75,0.4}
\definecolor{colorpppp}{rgb}{0.7,0.0,0.2}
\definecolor{coloroooo}{rgb}{0.45,0.0,0.0}
\definecolor{colorqqqq}{rgb}{0.1,0.7,0}

\definecolor{coloraaaa}{rgb}{0.6,0.6,0.6}

\newcommand{\eqnb}{\begin{equation}}
\newcommand{\eqne}{\end{equation}}

\title{Inertial manifolds for second-order parabolic equations\\in one space dimension}
\author{Igor Kukavica}
\date{}

\begin{document}
\maketitle
\begin{center}
Department of Mathematics\\
University of Southern California\\
Los Angeles, CA 90089\\
\texttt{kukavica@math.usc.edu}
\end{center}
\begin{abstract}
Inertial manifolds are positively invariant Lipschitz manifolds of finite dimension which attract all solutions with an exponential rate. We prove existence of inertial manifolds for a pointwise preparation of equations of the form
\[
u_t-u_{xx}+f(x,u,u_x)=0
\]
with periodic boundary conditions. In particular, the statement applies to pointwise preparations of the Burgers equation $u_t-u_{xx}+uu_x=f(x)$ and the modified Burgers equation $u_t-u_{xx}+uu_x-u=0$. The manifolds are parametrized by uniformly spaced nodal values and have the exponential tracking property.
\end{abstract}
\noindent\textbf{Keywords:} inertial manifold, global attractor, Burgers equation.

\medskip
\noindent\textbf{Mathematics Subject Classification:} 35B42, 35B41, 35K55, 37L25.
\newpage

\section{Introduction}\label{sec:introduction}
Inertial manifolds are positively invariant Lipschitz manifolds of finite dimension which attract all solutions with an exponential rate. They were introduced by Foias, Sell, and Temam in \cite{FST} and studied in numerous works; see, for instance,~\cite{CFNT,FSTi}. In addition to containing the global attractor, they may have the exponential tracking property, which means that every solution is approached exponentially by a solution on the manifold. We refer to \cite{FSTi} for this property.

We consider equations of the form
\begin{equation}\label{eq:introoriginal}
u_t-u_{xx}+f(x,u,u_x)=0
\end{equation}
with periodic boundary conditions. We assume that $f$ is infinitely smooth and periodic in the first variable. Particular examples are the Burgers equation
\[
u_t-u_{xx}+uu_x=f(x)
\]
and the modified Burgers equation
\[
u_t-u_{xx}+uu_x-u=0.
\]
In this paper, we prove that a pointwise preparation of~\eqref{eq:introoriginal} possesses inertial manifolds of every sufficiently large odd dimension. In particular, the result applies to a preparation of equations of the form
\[
u_t-u_{xx}+g^{(1)}(x,u)u_x+g^{(0)}(x,u)=0,
\]
but the dependence on $u_x$ in~\eqref{eq:introoriginal} need not be affine.

In general, the usual spectral gap condition does not hold for equations of the form~\eqref{eq:introoriginal}.

Our approach is related to~\cite{CFNT}, where the construction is based on the behavior of the Dirichlet quotient
\[
q(t)=\frac{\norm{v_x(\cdot,t)}_{L^2}^2}{\norm{v(\cdot,t)}_{L^2}^2},
\]
where $v$ is the difference of two solutions. For equations containing a first-order term, the estimate for $q'$ contains $q$ on the right-hand side. The gaps in the spectrum are then not sufficiently large to obtain the cone condition used in the spectral construction.

In~\cite{Kukavica2003}, we proved that global attractors for a class of dissipative equations in one space dimension are Lipschitz graphs above a finite number of Fourier modes. In~\cite{Kukavica2005}, the number of required modes was bounded algebraically in the parameters, improving the exponential bound in~\cite{Kukavica2003}. These results do not imply existence of an invariant manifold containing the attractor.

The main idea of our approach is to use a parametrization by uniformly spaced nodal values. The idea is related to determining nodes introduced by Foias and Temam in~\cite{FoiasTemam1984}; see also Foias and Titi~\cite{FTi}.

The preparation of the equation is chosen to depend pointwise on $u$ and~$u_x$. The difference $v=u_1-u_2$ of two solutions then satisfies
\[
v_t=v_{xx}+a(x,t)v_x+b(x,t)v
\]
with uniformly bounded coefficients. We may therefore use the zero-number property for scalar parabolic equations. Note that a cutoff depending on an integral norm would in general introduce nonlocal terms in the equation for~$v$.

In order to construct the manifold, we need two estimates. First, suppose that $v$ vanishes at $n$ prescribed nodes at time~$T$. By the zero-number theorem, $v$ has at least $n$ distinct zeros at every earlier positive time. A compactness argument shows that, if $n$ is sufficiently large, the $H^1$ norm decreases by a fixed factor on every unit time interval before~$T$. A second compactness argument bounds the initial $H^1$ norm by the initial high Fourier modes. We thus obtain an estimate which is uniform in~$T$.

Next, we consider differences of solutions whose initial values belong to the space of the first $n$ Fourier modes. We write each solution as its decaying linear part and a uniformly bounded remainder. While both linear parts are large, the preparation is active only on a small set. An elementary estimate for trigonometric polynomials, combined with the negative term in the Dirichlet-quotient identity, bounds the increase of the quotient on this interval. When one linear part is bounded and the other is large, the size of the difference gives an integrable bound. Once both are bounded, smoothing and the gauge argument of~\cite{Kukavica2003} apply. Compactness and the zero-number property then bound the inverse of the nodal map. The images of the Fourier space are consequently graphs with a uniform Lipschitz constant, and the terminal-node estimate shows that these graphs converge to the required manifold.

Vukadinovi\'c~\cite{Vukadinovic} proved existence of inertial manifolds for a class of diffusive Burgers equations, including the forced Burgers and Burgers-Sivashinsky equations. His use of the Cole-Hopf transform overcomes the spectral gap obstruction by exploiting the Burgers nonlinearity. This construction applies to that class of equations and does not treat a general nonlinearity $f(x,u,u_x)$ as considered here.

Kostianko and Zelik~\cite{KZ} established inertial-manifold results for scalar periodic reaction-diffusion-advection equations with affine dependence on~$u_x$. In Section~6.1, they extend their approach to general dependence on $u_x$ by adjoining the spatial derivative as an independent unknown. Their construction gives an inertial manifold for the resulting system in an enlarged phase space. To obtain an inertial manifold for the scalar equation, an additional argument is needed: the intersection with the derivative constraint need not be a manifold, and projection onto the scalar component need not be one-to-one. Thus, the construction in the enlarged space does not by itself settle the inertial-manifold problem for general gradient dependence in the original scalar phase space. In the present paper, we construct the manifold for the pointwise prepared equation directly in the scalar phase space, with nodal parametrization and exponential tracking.

The main result concerns the prepared equation. If the original equation has a bounded absorbing set in $C^1$, the preparation may be chosen to agree with it in a neighborhood of that set. The original attractor is then contained in the manifold, and the original solutions admit eventual exponential tracking, as stated in Section~\ref{sec:originalapplication}.

The organization of the paper is as follows. In Section~2, we prepare the equation and state the main result on existence of inertial manifolds. In Section~3, we derive the linear estimates needed in later sections and recall the zero-number property. In Section~4, we prove an estimate for solutions vanishing at the prescribed nodes at a terminal time. Section~5 shows the connection between the $H^1$ norm and the nodal values. Section~6 contains the construction of the inertial manifolds. In Section~7, we state the conclusions for the original equation.

\section{The main result}\label{sec:main}
We consider the equation
\begin{equation}\label{eq:original}
u_t-u_{xx}+f(x,u,u_x)=0
\end{equation}
with periodic boundary conditions with respect to $\Omega=[0,1]$. We assume that $f$ is infinitely smooth in all variables and periodic in $x$ with period one. We identify periodic functions on $\Omega$ with functions on $\T=\R/\mathbb Z$. Denote
\[
H^k=H^k_{\mathrm{per}}(\Omega)\qquad\text{and}\qquad H=H^0=L^2_{\mathrm{per}}(\Omega).
\]
The scalar product in $H^1$ is
\[
(u_1,u_2)_{H^1}=(u_1,u_2)+(u_{1x},u_{2x}),
\]
where $(\cdot,\cdot)$ denotes the scalar product in~$H$. Throughout the paper all functions are real valued.

We assume that the original equation~\eqref{eq:original} generates a global solution semigroup in the phase space under consideration and is dissipative in the following sense. There exists $R_{\mathrm d}>0$ such that, for every $R>0$, there exists $t_R>0$ for which every solution with admissible initial datum satisfies
\begin{equation}\label{eq:originaldissipation}
\norm{u(0)}_{L^2}\le R\quad\Longrightarrow\quad
\norm{u(t)}_{H^1}\le R_{\mathrm d},\qquad t\ge t_R.
\end{equation}
For the modified Burgers equation, the original phase space is restricted to initial data of zero spatial mean.

We prepare the equation as follows. Fix $R>0$ and choose a function $\chi\in C_c^\infty(\R^2)$ such that $0\le\chi\le1$ and $\chi=1$ on~$[-R,R]^2$. Denote
\begin{equation}\label{eq:cutoff}
F(x,u,p)=\chi(u,p)\bigl(u-f(x,u,p)\bigr),
\end{equation}
and consider the prepared equation
\begin{equation}\label{eq:prepared}
u_t-u_{xx}+u=F(x,u,u_x).
\end{equation}
Note that the dynamics of the equation~\eqref{eq:prepared} agrees with that of~\eqref{eq:original} as long as
\[
\norm{u(\cdot,t)}_{L^\infty}\le R\qquad\text{and}\qquad \norm{u_x(\cdot,t)}_{L^\infty}\le R.
\]

From here on, we consider the equation~\eqref{eq:prepared} with periodic boundary conditions. Assume that
\begin{enumerate}[label=(A\arabic*),leftmargin=*]
\item $F$ is infinitely smooth in $(x,u,p)$ and periodic in $x$ with period one;
\item $F$ and all its derivatives are bounded on $\T\times\R^2$;
\item There exists $R_*>0$ such that $F(x,u,p)=0$ whenever $|u|\ge R_*$.
\end{enumerate}
These assumptions hold for~\eqref{eq:cutoff}. Denote
\begin{equation}\label{eq:constants}
K=1+\norm F_{L^\infty}+\norm{F_u}_{L^\infty}+\norm{F_p}_{L^\infty}.
\end{equation}
The symbol $C$ denotes a positive constant whose value may change from one inequality to another.

Let $n\ge3$ be odd, and denote
\[
E_n=\spanop\bigl\{1,\cos(2\pi kx),\sin(2\pi kx):k=1,\ldots,(n-1)/2\bigr\}.
\]
Denote by $P_n$ the orthogonal projection onto $E_n$, and let $Q_n=I-P_n$. Denote
\[
x_j=\frac{j}{n},\quad j=0,\ldots,n-1,\qquad\text{and}\qquad \psi{}u=(u(x_0),\ldots,u(x_{n-1})).
\]
The restriction $\psi{}|_{E_n}$ is an isomorphism. Denote its inverse by $\theta_n=(\psi{}|_{E_n})^{-1}$. We use the usual Euclidean norm on~$\R^n$.

\begin{theorem}[inertial manifold for the prepared equation]\label{thm:main}
The equation~\eqref{eq:prepared} generates a global continuous solution semigroup $S(t)$ on~$H^1(\T)$. There exists a constant $n_0=n_0(K)$ such that for every odd $n\ge n_0$ there exists a mapping $\phi\colon \R^n\to H^1(\T)$ with the following properties:
\begin{enumerate}[label=(\arabic*),leftmargin=*]
\item The set $\calM=\phi(\R^n)$ is invariant under the flow, i.e., $S(t)\calM=\calM$ for all $t\ge0$.
\item The mapping $\phi$ is globally Lipschitz, and $\psi{}\phi$ is the identity on~$\R^n$. The set $\calM$ is a closed Lipschitz graph over $E_n$ with respect to the projection $\theta_n\psi{}$, whose kernel is~$\ker\psi{}$.
\item There exists a constant $L$ such that
\[
\norm{\phi(y)-\phi(z)}_{H^1}\le L|y-z|\qquad\text{and}\qquad |y-z|\le L\norm{\phi(y)-\phi(z)}_{L^2}.
\]
\item The global attractor $\calA$ of~\eqref{eq:prepared} is contained in~$\calM$.
\item For each $u_0\in H^1$, there exists $w_0\in\calM$ such that
\begin{equation}\label{eq:tracking}
\norm{S(t)u_0-S(t)w_0}_{H^1}
 \le C e^{-t/2}\bigl(\norm{Q_nu_0}_{H^1}+C_F\bigr),\qquad t\ge0.
\end{equation}
\end{enumerate}
\end{theorem}

Here $L$ depends on $F,n$, the constant $C$ depends only on $K,n$, and $C_F$ depends only on~$K$. All the constants are independent of the initial datum.

The property~(5) is called the exponential tracking property, or asymptotic completeness. 

Let $u_1$ and $u_2$ be two solutions of the prepared equation~\eqref{eq:prepared}, and denote their difference by $v=u_1-u_2$. Then
\begin{equation}\label{eq:linear}
v_t=v_{xx}+a(x,t)v_x+b(x,t)v,\qquad\text{with }\norm a_{L^\infty}+\norm b_{L^\infty}\le K,
\end{equation}
where, with $u_\theta=u_2+\theta(u_1-u_2)$,
\[
a=\int_0^1F_p(x,u_\theta,\partial_xu_\theta)\,d\theta\qquad\text{and}\qquad b=-1+\int_0^1F_u(x,u_\theta,\partial_xu_\theta)\,d\theta.
\]

\section{Linear estimates and the zero-number property}\label{sec:linear}

\subsection{Energy and smoothing estimates}
Let $U(t,s)$ denote the evolution operator for~\eqref{eq:linear}. The energy estimates and the representation formula for the heat equation give
\begin{equation}\label{eq:l2growth}
\norm{U(t,s)}_{L^2\to L^2}\le e^{C_K(t-s)}.
\end{equation}
For the $H^1$ norm on intervals of length at most one, we have
\begin{equation}\label{eq:h1growth}
\norm{U(t,s)}_{H^1\to H^1}\le M_K,\qquad0\le t-s\le1.
\end{equation}
The smoothing estimate at unit time reads
\begin{equation}\label{eq:smoothing}
\norm{U(s+1,s)}_{L^2\to H^1}\le M_K.
\end{equation}
Multiplying the equation by $v-v_{xx}$ and integrating leads to
\[
\begin{aligned}
 \frac12\frac d{dt}\norm v_{H^1}^2
 &=-\norm{v_x}_{L^2}^2-\norm{v_{xx}}_{L^2}^2+(av_x+bv,v-v_{xx})\\[4pt]
 &\le C_K\norm v_{H^1}^2-\frac12\norm{v_{xx}}_{L^2}^2.
\end{aligned}
\]
The $L^2$ estimate also gives a bound for $\int_s^{s+1/2}\norm{v_x}^2_{L^2}\,dt$. We choose a time in this interval at which the $H^1$ norm satisfies the corresponding bound. Applying~\eqref{eq:h1growth} from this time, we obtain~\eqref{eq:smoothing}.

We shall also use the following compactness property. For each $M>0$, solutions on $[0,1]$ with initial $H^1$ norm at most $M$ and coefficients bounded by $K$ are relatively compact in
\begin{equation}\label{eq:compactness}
C([\delta,1];C^1(\T))\quad\hbox{for every }\delta\in(0,1).
\end{equation}
Indeed, $av_x+bv$ is bounded in~$L^\infty_tL^2_x$. The representation formula for the heat equation gives uniform $H^s$ bounds for every $s<2$ on time intervals bounded away from zero, since $r^{-s/2}$ is integrable near zero. Choose $3/2<s'<s<2$. The equation gives a bound for the time derivative in~$H^{s-2}$. By interpolation, the solutions are equicontinuous in time with values in~$H^{s'}$. The compact embedding then gives~\eqref{eq:compactness}. Passing to a subsequence, we may assume that $a,b$ converge weak-*. Since $v,v_x$ converge strongly away from time zero, the limit satisfies the equation. The forcing is bounded on each such interval, so the usual local estimates for the heat equation give $W_p^{2,1}$ regularity for every finite~$p$. The initial values pass to the limit in~$L^2$. If the initial data converge strongly in $H^1$, the limit is the unique solution with that initial value.

Next, we note that a nonzero solution cannot vanish identically at a finite later time. The Dirichlet quotient for~\eqref{eq:linear} is
\[
q=\frac{\norm{v_x}_{L^2}^2}{\norm v_{L^2}^2}.
\]
On an interval where the denominator is nonzero, we obtain
\[
q'\le 2K^2(q+1)\qquad\text{and}\qquad \frac d{dt}\log\norm v_{L^2}\ge-q-K\sqrt q-K.
\]
Integrating these inequalities from a positive time, we conclude that a nonzero solution cannot vanish identically at a later time. By approximation, the same assertion holds for strong solutions with bounded coefficients.

\subsection{The zero-number property}
We use the scalar zero-number theorem of Angenent~\cite{Angenent}. Since the coefficient limits in the compactness argument are only bounded, we use the bounded lower-order coefficient formulation recalled in~\cite[Remark 1.3]{Lou}. For the unique-continuation arguments, we refer to Poon~\cite{Poon} and to the treatment using self-similar variables in~\cite{Kukavica1999}. For nontrivial strong solutions of
\[
|v_t-v_{xx}|\le K(|v_x|+|v|),
\]
the number of distinct spatial zeros is finite at interior times and is nonincreasing in time. It decreases strictly across a multiple zero. In particular, every interior time interval contains a time at which all zeros are simple. The result on an interval applies locally on the circle by choosing a point where the solution is nonzero and cutting the circle at that point. Both boundary values then remain nonzero for a short time.

Denote by $Z(v(t))$ the number of distinct zeros on the circle. Note that, if a function has finitely many simple zeros, this number is unchanged under sufficiently small $C^1$ perturbations. Indeed, outside disjoint neighborhoods of the zeros the function is bounded away from zero, and in each of these neighborhoods its derivative has a fixed nonzero sign.

If the initial datum $p$ belongs to $E_n$, the solution has at most $n-1$ zeros at every positive time. To include the case of multiple zeros at time zero, extend the solution to negative times by the finite Fourier sum $e^{t\partial_{xx}}p$, taking $a=b=0$ for $t<0$. This is a strong solution across time zero, and the zero-number theorem applies. 

\section{An estimate for solutions vanishing at prescribed nodes}\label{sec:B}
In order to construct the inertial manifold, we need the following lemma.
\begin{lemma}\label{lem:B}
There exists $n_0=n_0(K)$ such that, for every odd $n\ge n_0$, every solution of~\eqref{eq:linear} satisfies
\begin{equation}\label{eq:B}
\norm{v(t)}_{H^1}\le C_{K,n}e^{-t/2}\norm{Q_nv(0)}_{H^1},\qquad0\le t\le T,
\end{equation}
whenever $\psi{}v(T)=0$ and $T\ge1$.
\end{lemma}
\begin{proof}[Proof of Lemma~\ref{lem:B}]
We first show that, for each fixed odd $n$, there exists $D_{K,n}$ such that
\begin{equation}\label{eq:initialangle}
\norm{v(0)}_{H^1}\le D_{K,n}\norm{Q_nv(0)}_{H^1},
\end{equation}
whenever $\psi{}v(T)=0$ and $T\ge1$.
Assume that $v$ is nonzero and that $\psi{}v(T)=0$. By the zero-number theorem, we have
\begin{equation}\label{eq:manyzeros}
Z(v(t))\ge n\quad(0<t<T).
\end{equation}

Assume, contrary to the assertion, that~\eqref{eq:initialangle} fails. By linearity, there is a sequence of solutions such that
\[
\norm{v_\ell(0)}_{H^1}=1\qquad\text{and}\qquad\norm{Q_nv_\ell(0)}_{H^1}\longrightarrow0,\qquad\text{and}\qquad T_\ell\ge1.
\]
Since $E_n$ is finite dimensional and $\norm{P_nv_\ell(0)}_{H^1}\le1$, we may pass to a subsequence for which $P_nv_\ell(0)\to p$ strongly in $H^1$, with $p\in E_n$. Together with $Q_nv_\ell(0)\to0$ in $H^1$, this gives $v_\ell(0)\to p$ strongly in $H^1$ and $\norm p_{H^1}=1$. By the compactness property on $[0,1]$, we obtain a solution $w$ with $w(0)=p$ and convergence in $C^1$ at every time in~$(0,1)$. Extend $w$ backward by $e^{t\partial_{xx}}p$ for negative times. For these negative times, the solution has at most $n-1$ zeros. Choose a sufficiently small $t_*\in(0,1)$ at which $w(t_*)$ has only simple zeros. There are at most $n-1$ such zeros. Since simple zeros are preserved under small $C^1$ perturbations, we obtain $Z(v_\ell(t_*))\le n-1$ for large $\ell$, contradicting~\eqref{eq:manyzeros}. This proves~\eqref{eq:initialangle}.

Next, we show that there exists $N_K$ such that
\begin{equation}\label{eq:unitcontraction}
\norm{v(1)}_{H^1}\le\tfrac12\norm{v(0)}_{H^1},
\end{equation}
whenever $Z(v(t))\ge N_K$ for every $0<t<1$.
Assume, contrary to the assertion, that there is a sequence of solutions $\{w_\ell\}_{\ell=1}^\infty$ with $\norm{w_\ell(0)}_{H^1}=1$, at least $\ell$ zeros throughout $(0,1)$, and $\norm{w_\ell(1)}_{H^1}>1/2$. By~\eqref{eq:compactness} and a diagonal subsequence argument, we obtain a limit $w$ such that
\[
w_\ell\longrightarrow w\quad\text{in }C([\delta,1];C^1(\T))\quad\text{for every }0<\delta<1.
\]
In particular, the convergence at time one is strong in $H^1$, so $\norm{w(1)}_{H^1}\ge1/2$. By the zero-number theorem for bounded coefficients, there is a time $t_*\in(0,1)$, sufficiently close to one, at which $w$ has finitely many simple zeros. The convergence in $C^1$ implies that $w_\ell(t_*)$ has the same number of zeros for all sufficiently large~$\ell$. This is a contradiction.

Choose an odd $n\ge3$ such that $n\ge N_K$. By~\eqref{eq:manyzeros}, we may apply~\eqref{eq:unitcontraction} on every unit interval contained in $[0,T]$. Denote $k=\lfloor t\rfloor$. Applying~\eqref{eq:h1growth} on the remaining interval, we get
\[
\norm{v(t)}_{H^1}
 \le M_K2^{-k}\norm{v(0)}_{H^1}
 \le2M_Ke^{-t/2}\norm{v(0)}_{H^1}.
\]
The assertion now follows from~\eqref{eq:initialangle}. The case $v=0$ is immediate.
\end{proof}

\begin{remark}
Replacing $1/2$ in~\eqref{eq:unitcontraction} by any prescribed $\rho\in(0,1)$ gives the rate $-\log\rho$ in~\eqref{eq:B} after increasing~$n$.
\end{remark}

\section{Nodal values and the \texorpdfstring{$H^1$}{H1} norm}\label{sec:A}
We next estimate the difference of two solutions with initial data in $E_n$ by its $L^2$ norm and by its values at the prescribed nodes.
\begin{lemma}\label{lem:A}
Let $n\ge3$ be odd, and let $u_i(t)=S(t)p_i$, where $p_i\in E_n$ and $i=1,2$. Denote $v=u_1-u_2$. Then
\begin{equation}\label{eq:Alinear}
\norm{v(t)}_{H^1}\le L_n\norm{v(t)}_{L^2}\qquad\text{and}\qquad
\norm{v(t)}_{H^1}\le L_n|\psi{}v(t)|,\qquad t\ge1.
\end{equation}
\end{lemma}
The constant $L_n$ depends on the fixed nonlinearity $F$ and on $n$, and is independent of $p_1,p_2$, and~$t$. We first recall a standard sublevel estimate for trigonometric polynomials. It follows from the Tur\'an-Nazarov inequality; see~\cite[Theorem~I]{Nazarov}. We include an elementary proof for completeness.

\begin{lemma}\label{lem:sublevel}
For each odd $n\ge3$, there exists $C_n$ such that
\begin{equation}\label{eq:sublevel}
\bigl|\{x\in\T:|p(x)|\le d\}\bigr|
\le C_n\left(\frac d{\norm p_{L^2}}\right)^{1/n},
\qquad p\in E_n\setminus\{0\},\quad d>0.
\end{equation}
Here $|D|$ denotes the Lebesgue measure of $D\subseteq\T$.
\end{lemma}
\begin{proof}[Proof of Lemma~\ref{lem:sublevel}]
Without loss of generality, $\norm p_{L^2}=1$. The mapping from $E_n$ to $\R^{n}$ which assigns to $p$ its derivatives of orders zero through $n-1$ at a fixed point is an isomorphism. Consequently, at every point, at least one of these derivatives has absolute value at least $c_n>0$. All derivatives of normalized elements of $E_n$ are bounded uniformly. Thus, $\T$ can be covered by a bounded number of intervals on each of which one derivative of order $k$, $0\le k\le n-1$, has absolute value at least~$c_n/2$.

For $k=0$, the sublevel set on this interval is empty when $d<c_n/2$. For $k\ge1$, let $s$ be its measure. Choose $k+1$ points of the sublevel set whose mutual distances are at least a fixed multiple of~$s/k$. The divided-difference formula and the generalized mean-value theorem give
\[
\frac{c_n}{2}\le C_k d s^{-k}.
\]
Hence, $s\le C_n d^{1/k}$. Since $k\le n-1$, we have $d^{1/k}\le d^{1/n}$ for $0<d<1$. Summing over the covering intervals proves~\eqref{eq:sublevel} when $d$ is small. Increasing $C_n$ covers the remaining values of~$d$.
\end{proof}

\begin{lemma}\label{lem:quotient}
For the solutions in Lemma~\ref{lem:A}, we have
\begin{equation}\label{eq:uniformquotient}
q(t)=\frac{\norm{v_x(t)}_{L^2}^2}{\norm{v(t)}_{L^2}^2}\le C_{F,n},\qquad t\ge0,
\end{equation}
provided $p_1\ne p_2$.
\end{lemma}
\begin{proof}[Proof of Lemma~\ref{lem:quotient}]
Choose $R_*>0$ such that $F(x,u,p)=0$ whenever $|u|\ge R_*$. Such a number exists by~\eqref{eq:cutoff}. The prepared equation itself satisfies the dissipative estimate
\begin{equation}\label{eq:prepareddissipation}
\frac d{dt}\norm u_{H^1}^2+\norm u_{H^1}^2\le K^2,
\end{equation}
which we obtain by multiplying the prepared equation~\eqref{eq:prepared} by $u-u_{xx}$, integrating, and using
\[
\frac12\frac d{dt}\norm u_{H^1}^2+\norm{u-u_{xx}}_{L^2}^2
\le\frac12\norm{u-u_{xx}}_{L^2}^2+\frac12 K^2.
\]
Here
\[
\norm{u-u_{xx}}_{L^2}^2
=\norm u_{L^2}^2+2\norm{u_x}_{L^2}^2+\norm{u_{xx}}_{L^2}^2
\ge\norm u_{H^1}^2.
\]
In particular, above the radius $\sqrt2 K$, the $H^1$ norm decreases exponentially. A fixed annulus above this radius is crossed in a bounded time, and a ball of radius greater than $K$ is forward invariant. Note that these are consequences of the preparation.

Write
\[
u_i(t)=z_i(t)+r_i(t),\qquad\text{where}\qquad z_i(t)=e^{-t(1-\partial_{xx})}p_i.
\]
The representation formula and the bound on $F$ give
\begin{equation}\label{eq:quotientremainder}
\norm{r_i(t)}_{L^\infty}\le K\qquad\text{and}\qquad
\norm{r_i(t)}_{H^1}\le C K,\qquad t\ge0.
\end{equation}
Moreover, $z_i(t)\in E_n$ and
\begin{equation}\label{eq:linearbackward}
\norm{z_i(s)}_{L^2}\ge e^{t-s}\norm{z_i(t)}_{L^2},\qquad0\le s\le t.
\end{equation}

Denote
\[
\mathcal L=-\partial_{xx}\qquad\text{and}\qquad
G=F(x,u_1,(u_1)_x)-F(x,u_2,(u_2)_x).
\]
Set
\[
w=\frac{v}{\norm v_{L^2}}
\qquad\text{and}\qquad
D=\norm{(\mathcal L-q)w}_{L^2}.
\]
Differentiating the quotient gives
\begin{equation}\label{eq:quotientidentity}
q'=-2D^2+2\bigl((\mathcal L-q)w,G/\norm v_{L^2}\bigr).
\end{equation}
Consequently,
\begin{equation}\label{eq:quotientbounds}
q'\le-D^2+\frac{\norm G_{L^2}^2}{\norm v_{L^2}^2}
\qquad\text{and}\qquad
q'\le\frac{\norm G_{L^2}^2}{2\norm v_{L^2}^2}.
\end{equation}
The global Lipschitz bound for $F$ also implies $q'\le C_K(1+q)$.

First, we consider the interval on which both linear parts are large. Define
\[
D_i(t)=\{x\in\T:|u_i(x,t)|\le R_*\}\qquad\text{and}\qquad
\omega(t)=|D_1(t)\cup D_2(t)|.
\]
The function $G$ vanishes outside this union, and
\[
|G|\le K(|v|+|v_x|).
\]
Since $\norm w_{L^2}=1$, the one-dimensional interpolation inequalities yield
\[
\norm w_{L^\infty}^2\le1+2\sqrt q\qquad\text{and}\qquad
\norm{w_x}_{L^\infty}^2\le q+2\sqrt q\sqrt{D^2+q^2}.
\]
Here we used $\norm{w_{xx}}_{L^2}^2=D^2+q^2$. If $q\le Q$, the first inequality in~\eqref{eq:quotientbounds} thus gives
\begin{equation}\label{eq:localizedquotient}
q'\le-D^2+C_{K,Q}\omega(t)(1+D)\le C_{K,Q}\omega(t).
\end{equation}
The last step follows by completing the square and using $0\le\omega\le1$.

By~\eqref{eq:quotientremainder}, $D_i(t)$ is contained in the set where $|z_i(t)|\le R_*+K$. Lemma~\ref{lem:sublevel} gives
\begin{equation}\label{eq:activeset}
\omega(t)\le C_n\sum_{i=1}^2
\left(\frac{R_*+K}{\norm{z_i(t)}_{L^2}}\right)^\alpha,
\qquad\text{where}\qquad \alpha=\frac1n,
\end{equation}
whenever both denominators are nonzero.

Choose a number $L>R_*+K+1$, to be fixed below. If both initial linear parts have norm greater than $L$, let $t_1$ be the first time that their smaller norm is $L$, and otherwise set $t_1=0$. In the first case,~\eqref{eq:linearbackward} gives
\[
\norm{z_i(t)}_{L^2}\ge Le^{t_1-t},\qquad0\le t\le t_1,\quad i=1,2.
\]
Thus,
\begin{equation}\label{eq:activeintegral}
\int_0^{t_1}\omega(t)\,dt
\le\frac{C_n}{\alpha}\left(\frac{R_*+K}{L}\right)^\alpha.
\end{equation}
Initially, $q(0)\le\pi^2(n-1)^2$. Fix $Q=\pi^2(n-1)^2+2$, and then choose $L=L(F,n)$ so large that the product of the right side of~\eqref{eq:activeintegral} and the constant in~\eqref{eq:localizedquotient} is less than one. A first-crossing argument shows that
\[
q(t)\le\pi^2(n-1)^2+1,\qquad0\le t\le t_1.
\]
The same conclusion at $t_1=0$ is immediate.

Next, we consider the mixed interval. Relabel the solutions so that $\norm{z_1(t_1)}_{L^2}\le L$, and put $H=L+2K+1$. For all $t\ge t_1$, we have $\norm{u_1(t)}_{L^2}\le L+K$. If $\norm{z_2(t_1)}_{L^2}>2H$, let $t_2$ be the first time at which $\norm{z_2(t_2)}_{L^2}=2H$; otherwise put $t_2=t_1$. On $[t_1,t_2]$ in the nonempty case,
\[
\norm{v(t)}_{L^2}\ge\norm{z_2(t)}_{L^2}-H\ge\tfrac12\norm{z_2(t)}_{L^2}.
\]
Since $\norm G_{L^2}\le2K$, the second inequality in~\eqref{eq:quotientbounds} and~\eqref{eq:linearbackward} imply
\begin{equation}\label{eq:mixedquotient}
q(t)-q(t_1)\le8K^2\int_{t_1}^{t}\frac{ds}{\norm{z_2(s)}_{L^2}^2}
\le\frac{K^2}{H^2},\qquad t_1\le t\le t_2.
\end{equation}

At time $t_2$, both linear parts are bounded in~$L^2$. Since they belong to $E_n$,~\eqref{eq:quotientremainder} bounds both solutions in $H^1$ at $t_2$ and at every later time. Parabolic regularization gives uniform bounds for their spatial derivatives through order three at times $t\ge t_2+1$. Their equations then give bounds for $\partial_tu_i$ and~$\partial_x\partial_tu_i$. All these bounds depend only on~$F,n$. On $[t_2,t_2+1]$, the inequality $q'\le C_K(1+q)$ bounds the quotient uniformly.

For completeness, we recall the argument of~\cite[Theorem~2.4]{Kukavica2003} for the remaining interval. Denote $s_0=t_2+1$. The coefficients $a,b$ in~\eqref{eq:linear} satisfy
\[
\sup_{t\ge s_0}\bigl(\norm a_{L^\infty}+\norm b_{L^\infty}
+\norm{a_x}_{L^\infty}+\norm{a_t}_{L^\infty}\bigr)\le C_{F,n}.
\]
Set
\[
\bar a(t)=\int_0^1a(x,t)\,dx\qquad\text{and}\qquad
h(x,t)=\tfrac12\int_0^x(\bar a(t)-a(y,t))\,dy.
\]
The periodic function $W=e^{-h}v$ solves
\begin{equation}\label{eq:gaugeequation}
W_t=W_{xx}+\bar a(t)W_x+c(x,t)W,
\qquad\text{where}\qquad c=h_{xx}+h_x^2+ah_x+b-h_t.
\end{equation}
The functions $h,h_x,c$ are uniformly bounded. Denote
\[
Q_W(t)=\frac{\norm{W_x(t)}_{L^2}^2}{\norm{W(t)}_{L^2}^2}.
\]
Periodicity gives $(W,W_x)=0$ and $(\mathcal LW,W_x)=0$. Thus, the spatially constant drift contributes zero when we differentiate~$Q_W$. We obtain
\[
\begin{aligned}
Q_W'
&=-2\frac{\norm{(\mathcal L-Q_W)W}_{L^2}^2}{\norm W_{L^2}^2}
+2\frac{\bigl((\mathcal L-Q_W)W,cW\bigr)}{\norm W_{L^2}^2}\\
&\le-\frac{\norm{(\mathcal L-Q_W)W}_{L^2}^2}{\norm W_{L^2}^2}
+\norm c_{L^\infty}^2.
\end{aligned}
\]
The bounds on $q(s_0),h,h_x$ also bound $Q_W(s_0)$ uniformly. Let $\lambda_j=(2\pi j)^2$ for $j\ge0$. Choose $j$ so large that
\[
\Lambda=\frac{\lambda_j+\lambda_{j+1}}2>Q_W(s_0)
\qquad\text{and}\qquad
\frac{\lambda_{j+1}-\lambda_j}2>\sup_{t\ge s_0}\norm{c(t)}_{L^\infty}.
\]
Write $W(x)=\sum_{k\in\mathbb Z}\widehat W(k)e^{2\pi ikx}$. Since every eigenvalue of $\mathcal L$ lies at distance at least $(\lambda_{j+1}-\lambda_j)/2$ from $\Lambda$, Parseval's identity gives
\[
\begin{aligned}
\norm{(\mathcal L-\Lambda)W}_{L^2}^2
&=\sum_{k\in\mathbb Z}\bigl((2\pi k)^2-\Lambda\bigr)^2|\widehat W(k)|^2\\
&\ge\left(\frac{\lambda_{j+1}-\lambda_j}{2}\right)^2\norm W_{L^2}^2.
\end{aligned}
\]
Consequently, at $Q_W=\Lambda$ we have
\[
Q_W'\le-\left(\frac{\lambda_{j+1}-\lambda_j}2\right)^2
+\sup_{t\ge s_0}\norm{c(t)}_{L^\infty}^2<0.
\]
Hence, $Q_W(t)\le\Lambda$ for $t\ge s_0$. The bounds on $h,h_x$ return this estimate to $v$, since
\[
q(t)\le2e^{4\|h(t)\|_{L^\infty}}\bigl(Q_W(t)+\|h_x(t)\|_{L^\infty}^2\bigr).
\]
Together with~\eqref{eq:localizedquotient}--\eqref{eq:mixedquotient}, this proves~\eqref{eq:uniformquotient}.
\end{proof}

\begin{proof}[Proof of Lemma~\ref{lem:A}]
The first inequality follows from Lemma~\ref{lem:quotient}, in fact for every $t\ge0$. Multiplying~\eqref{eq:linear} by $v$ and using the quotient bound gives
\begin{equation}\label{eq:lognormgrowth}
\left|\frac d{dt}\log\norm{v(t)}_{L^2}\right|\le\beta_{F,n}.
\end{equation}
Suppose that the second inequality in~\eqref{eq:Alinear} fails. Then there are pairs of solutions as above, times $T_j\ge1$, and normalized differences
\[
w_j(s)=\frac{v_j(T_j+s)}{\norm{v_j(T_j)}_{L^2}},\qquad-\tfrac12\le s\le\tfrac12,
\]
such that
\begin{equation}\label{eq:nodalnormalize}
\norm{w_j(0)}_{L^2}=1\qquad\text{and}\qquad |\psi{}w_j(0)|\longrightarrow0.
\end{equation}
By~\eqref{eq:uniformquotient} and~\eqref{eq:lognormgrowth}, the $H^1$ norms of $w_j$ are uniformly bounded on this interval. The compactness property~\eqref{eq:compactness} gives, after passing to a subsequence,
\[
w_j\longrightarrow w\quad\text{in }C([-\tfrac14,\tfrac14];C^1(\T)).
\]
The limit solves a scalar parabolic equation with bounded coefficients and satisfies
\[
\norm{w(0)}_{L^2}=1\qquad\text{and}\qquad w(0,x_k)=0,\quad0\le k<n.
\]
Choose $s_*\in(-1/4,0)$ at which all zeros of $w(s_*)$ are simple. The zero-number theorem implies that there are at least $n$ of them. Uniform convergence then gives at least $n$ zeros of $w_j(s_*)$ for large~$j$. On the other hand, the initial differences belong to $E_n$, so every $w_j(s_*)$ has at most $n-1$ zeros. This contradiction proves the nodal estimate. For each fixed $t>0$, the same argument on a time window of radius $\min\{t/2,1/2\}$ gives the estimate with a constant which may depend on~$t$. This also proves the Lipschitz assertion for every positive time used below.
\end{proof}

\section{Construction of the inertial manifold for the prepared equation}\label{sec:construction}
We first note that the equation~\eqref{eq:prepared} is globally well-posed in~$H^1$. Indeed, the mapping $u\mapsto F(x,u,u_x)$ is globally Lipschitz from $H^1$ to $L^2$ and is bounded in~$L^2$. The representation formula for the heat equation, together with integrability of $r^{-1/2}$ near zero, gives a unique global mild solution continuous in~$H^1$. It also gives continuous dependence on the initial datum and regularity at positive times.

Denote $A=1-\partial_{xx}$. The representation formula
\begin{equation}\label{eq:duhamel}
S(t)u_0=e^{-tA}u_0+R_t(u_0)\quad\text{where}\quad R_t(u_0)=\int_0^te^{-(t-s)A}F(x,u(s),u_x(s))\,ds,
\end{equation}
gives
\begin{equation}\label{eq:remainder}
\sup_{t\ge0,u_0}\norm{R_t(u_0)}_{H^\sigma}\le C_{F,\sigma}
 \quad(0\le\sigma<2)\quad\text{and}\quad \norm{R_t(u_0)}_{L^\infty}\le K.
\end{equation}
The first estimate follows by integrating $Ce^{-r}(1+r^{-\sigma/2})K$ over $(0,\infty)$. Denote $C_F=C_{F,1}$.

\begin{proposition}\label{prop:graphs}
For every $t>0$, the set $\mathcal M_t=S(t)E_n$ is a graph parametrized by the nodal values, i.e.,
\[
\mathcal M_t=\phi_t(\R^n)\qquad\text{and}\qquad \psi{}\phi_t=\mathrm{id}.
\]
The mappings $\phi_t\colon \R^n\to H^1$ are globally Lipschitz. Moreover,
\begin{equation}\label{eq:tubes}
\norm{Q_n\phi_t(y)}_{H^1}\le C_F\quad\text{and}\quad \norm{\phi_t(y)-\theta_ny}_{H^1}\le D_n,\qquad t>0,\ y\in\R^n.
\end{equation}
\end{proposition}
For $t\ge1$, the Lipschitz constants are bounded uniformly in~$t$. The constants in~\eqref{eq:tubes} are independent of $t>0$ and $y\in\R^n$.
\begin{proof}[Proof of Proposition~\ref{prop:graphs}]
Denote $G_t=\psi{}S(t)|_{E_n}$. We first show that $G_t$ is one-to-one. If $G_t(p)=G_t(q)$, the difference of the solutions starting at $p,q$ vanishes at the $n$ prescribed nodes at time~$t$. Since the initial difference belongs to $E_n$, the zero-number theorem implies $p=q$. The mapping $G_t$ is also continuous. By~\eqref{eq:duhamel}, we have
\[
G_t(p)=\psi{}e^{-tA}p+\text{a uniformly bounded term}.
\]
The linear mapping $\psi{}e^{-tA}|_{E_n}$ is invertible. Consequently, $G_t$ is proper, i.e., the inverse image of every compact set is compact. By invariance of domain, its image is open. Properness implies that the image is closed. It follows that $G_t$ is onto $\R^n$ and is a homeomorphism. Denote $\phi_t=S(t)G_t^{-1}$.

Let two points of $\mathcal M_t$ be given. Their difference is covered by Lemma~\ref{lem:A}. Thus,
\begin{equation}\label{eq:graphlipschitz}
\norm{\phi_t(y)-\phi_t(z)}_{H^1}\le L_n|y-z|,\qquad t\ge1.
\end{equation}
The first estimate in that lemma also gives
\begin{equation}\label{eq:secantnorm}
\norm{\phi_t(y)-\phi_t(z)}_{H^1}
 \le L_n\norm{\phi_t(y)-\phi_t(z)}_{L^2}.
\end{equation}
Since $e^{-tA}p\in E_n$, the first bound in~\eqref{eq:tubes} follows from~\eqref{eq:remainder}. Applying $I-\theta_n\psi{}$ to~\eqref{eq:duhamel} gives the identity
\[
\phi_t(y)=\theta_ny+(I-\theta_n\psi{})R_t(G_t^{-1}y),
\]
which proves the second bound. The operator $\theta_n\psi{}$ is bounded on~$H^1$.
\end{proof}

\begin{proposition}\label{prop:limit}
For every sufficiently large odd $n$, there exists a globally Lipschitz map $\phi\colon \R^n\to H^1$ such that $\psi{}\phi=\mathrm{id}$ and
\[
\sup_{y\in\R^n}\norm{\phi_T(y)-\phi(y)}_{H^1}\le CC_Fe^{-T/2},\qquad T\ge1.
\]
Its image satisfies the graph and Lipschitz assertions of Theorem~\ref{thm:main}.
\end{proposition}
\begin{proof}[Proof of Proposition~\ref{prop:limit}]
Fix an odd $n$ sufficiently large for Lemma~\ref{lem:B}. For $T\ge1$ and $s\ge0$ and $y\in\R^n$, choose $p,q\in E_n$ with
\[
S(T)p=\phi_T(y)\qquad\text{and}\qquad S(T+s)q=\phi_{T+s}(y).
\]
Consider the solutions on $[0,T]$ with initial data $p$ and~$S(s)q$. Their nodal values agree at time $T$, and
\[
\norm{Q_n(p-S(s)q)}_{H^1}\le C_F.
\]
Lemma~\ref{lem:B} yields
\begin{equation}\label{eq:cauchy}
\sup_{y\in\R^n}\norm{\phi_{T+s}(y)-\phi_T(y)}_{H^1}
 \le C C_Fe^{-T/2}.
\end{equation}
It follows that the mappings converge uniformly in $H^1$ to a mapping $\phi\colon \R^n\to H^1$. Passing to the limit, we obtain $\psi{}\phi=\mathrm{id}$, the estimate~\eqref{eq:graphlipschitz} for $\phi$, and
\begin{align}
&\phi(y)=\theta_ny+r(y),\qquad\text{where }\psi{}r(y)=0,\notag\\
\shortintertext{with}
&\sup_y\norm{r(y)}_{H^1}\le D_n\qquad\text{and}\qquad\sup_y\norm{Q_n\phi(y)}_{H^1}\le C_F.\label{eq:limittube}
\end{align}
Thus, $\calM=\phi(\R^n)$ is a closed Lipschitz graph in~$H^1$. Passing to the limit in~\eqref{eq:secantnorm} and then applying the bounded map $\psi{}\colon H^1\to\R^n$ proves the $L^2$ Lipschitz bound for its inverse asserted in Theorem~\ref{thm:main}.
\end{proof}

\subsection{Invariance and exponential tracking}
We now prove the invariance and exponential tracking assertions of Theorem~\ref{thm:main}.
\begin{proof}[Completion of the proof of Theorem~\ref{thm:main}]
We first show that $\calM$ is positively invariant.
Let $\tau\ge0$. By the definition of $\phi_T$, we have
\[
S(\tau)\phi_T(y)=\phi_{T+\tau}\bigl(\psi{}S(\tau)\phi_T(y)\bigr).
\]
Let $T\to\infty$. Continuity of $S(\tau)$ and $\psi{}$, together with uniform convergence~\eqref{eq:cauchy} and the common Lipschitz constant, gives
\[
S(\tau)\phi(y)=\phi\bigl(\psi{}S(\tau)\phi(y)\bigr).
\]
Hence, $S(\tau)\calM\subseteq\calM$.

Next, we show that the restriction of $S(\tau)$ to $\calM$ is onto.
Denote $H_\tau(y)=\psi{}S(\tau)\phi(y)$. This mapping is continuous and one-to-one. Indeed, equality of the nodal values implies equality of the two images in $\calM$, and backward uniqueness then implies equality of the initial points. By~\eqref{eq:duhamel} and~\eqref{eq:limittube}, we have
\[
H_\tau(y)=\psi{}e^{-\tau A}\theta_ny+\text{a term bounded uniformly in }y.
\]
The displayed linear map is invertible. Thus, $H_\tau$ is proper and, by invariance of domain, onto. This proves $S(\tau)\calM=\calM$.

Now, we show existence of the global attractor $\calA$ and that $\calA\subseteq\calM$. The representation formula gives
\[
\norm{S(t)u_0}_{H^1}\le e^{-t}\norm{u_0}_{H^1}+C_F.
\]
Thus, the equation has an absorbing ball in~$H^1$. Taking $1<\sigma<2$ in~\eqref{eq:remainder} and using the smoothing estimate for the linear term, we obtain asymptotic compactness in~$H^1$. Consequently, the solution semigroup has a compact global attractor.

Let $u(t)$ be a solution on the attractor, defined and bounded for all $t\in\R$. Choose $p_R\in E_n$ such that
\[
S(R)p_R=\phi_R(\psi{}u(0)).
\]
We compare $u$ on $[-R,0]$ with the solution starting at $p_R$ at time~$-R$. At time zero their nodal values agree. By Lemma~\ref{lem:B}, we obtain
\[
\norm{u(0)-\phi_R(\psi{}u(0))}_{H^1}
 \le Ce^{-R/2}\norm{Q_nu(-R)}_{H^1}.
\]
Letting $R\to\infty$ proves $u(0)=\phi(\psi{}u(0))$. Hence, $\calA\subseteq\calM$.

Next, we establish exponential convergence to a solution on the manifold~$\calM$.
Let $u_0\in H^1$ be arbitrary. Since $S(T)\calM=\calM$ and $\calM$ is parametrized by the nodal values, for every $T\ge1$ there exists $w_T\in\calM$ such that
\[
\psi{}S(T)w_T=\psi{}S(T)u_0.
\]
Lemma~\ref{lem:B} and~\eqref{eq:limittube} imply, simultaneously for $0\le t\le T$,
\begin{equation}\label{eq:finite_tracking}
\norm{S(t)u_0-S(t)w_T}_{H^1}
 \le Ce^{-t/2}\bigl(\norm{Q_nu_0}_{H^1}+C_F\bigr).
\end{equation}
Taking $t=0$, we conclude that $w_T$ is bounded in~$H^1$. The nodal values are therefore bounded in~$\R^n$. Passing to a subsequence, we may assume that $T_j\to\infty$ and $w_{T_j}\to w_0\in\calM$ in~$H^1$. For each fixed $t\ge0$, we pass to the limit in~\eqref{eq:finite_tracking} along this sequence and obtain~\eqref{eq:tracking}. This completes the proof of Theorem~\ref{thm:main}.
\end{proof}

\section{The original equation}\label{sec:originalapplication}
For the original equation~\eqref{eq:original}, we use the following additional bound on $u$ and~$u_x$.

\begin{corollary}\label{cor:original}
Assume that the equation~\eqref{eq:original} generates a global solution semigroup in the phase space under consideration. Assume also that there exists $R_0>0$ such that every solution satisfies
\begin{equation}\label{eq:absorbingjet}
\norm{u(t)}_{L^\infty}+\norm{u_x(t)}_{L^\infty}\le R_0,
\qquad t\ge t_0(u_0).
\end{equation}
Choose $R_1>R_0$ and let $\chi=1$ on $[-R_1,R_1]^2$ in~\eqref{eq:cutoff}. Then, for every such solution, there exists a solution of the prepared equation on $\calM$ which approaches it exponentially for $t\ge t_0$. This latter solution satisfies the original equation for all sufficiently large times. Moreover, every bounded solution of the original equation defined for all real times, with $u,u_x\in[-R_0,R_0]$, is contained in~$\calM$. In particular, this holds for a global attractor of the original equation satisfying this bound.
\end{corollary}
\begin{proof}[Proof of Corollary~\ref{cor:original}]
For $t\ge t_0$, the condition~\eqref{eq:absorbingjet} implies that $u$ solves the prepared equation. Apply Theorem~\ref{thm:main} at time~$t_0$. Then there exists a solution $w(t)\in\calM$ of the prepared equation, defined for $t\ge t_0$, such that
\[
\norm{u(t)-w(t)}_{H^1}\le C(u_0)e^{-(t-t_0)/2}.
\]
The difference $u-w$ satisfies~\eqref{eq:linear}. Applying the smoothing estimate over a unit time interval, as in Section~\ref{sec:linear} with $3/2<s<2$, we obtain
\[
\norm{u(t)-w(t)}_{C^1}\le C_K\norm{u(t-1)-w(t-1)}_{H^1},\qquad t\ge t_0+1.
\]
It follows that $\norm{u(t)-w(t)}_{C^1}<R_1-R_0$ for all sufficiently large~$t$. Thus, both $w$ and $w_x$ belong to $[-R_1,R_1]$. Since $\chi=1$ there, $w$ satisfies the original equation for all such times.

A bounded solution of the original equation defined for all real times and satisfying the stated bounds on $u,u_x$ is also a bounded solution of the prepared equation for all real times. The argument in Section~\ref{sec:construction} showing $\calA\subseteq\calM$ then implies that this solution belongs to~$\calM$.
\end{proof}

\section*{Acknowledgment}
The author was supported in part by the Simons grant SFI-MPS-TSM-00025796.

The present paper is based on an unpublished manuscript written by the author in 2004, using nodal values to parametrize the inertial manifold. In revisiting and completing that manuscript, the author used ChatGPT 6.1 (OpenAI) to assist with drafting and polishing the paper.

\end{document}